\documentclass[12pt,letterpaper]{article}
\usepackage[utf8]{inputenc}
\usepackage{amsmath}
\usepackage{bm}
\usepackage{amsfonts}
\usepackage{amssymb}
\usepackage{tabu}
\usepackage{makecell}
\usepackage{wrapfig}
\usepackage[svgnames]{xcolor}
\usepackage[margin=0.75in]{geometry}
\usepackage{lipsum,booktabs}

\begin{document}
\begin{titlepage}
\title{Imaginary Number Bases}
\author{Philip Herd}
\maketitle

\begin{abstract}
An expansion upon Donald Kunth's quarter-imaginary base system is introduced to handle any imaginary number base where its real part is zero and the absolute value of its imaginary part is greater than one. A brief overview on number bases is given as well as conversion to both positive and negative bases. Additionally gives examples for addition, subtraction, multiplication for imaginary bases and adds a division method for imaginary bases as well as mentions possible uses.
\end{abstract}

\end{titlepage}

\section*{Introduction}
A number base (or radix) is a positional number system that uses a set of digits to represent numbers. This takes the form of $(\ldots d_{4}d_{3}d_{2}d_{1}d_{0}.d_{-1}d_{-2}d_{-3}\ldots)_{b}$ where the symbol between $d_{1}$ and $d_{-1}$ is the radix point, used to denote the whole and fractional parts of a number. This may also be either a comma or a raised dot depending on culture~\cite{knu2}. The most familiar base is base 10 (decimal) which uses ten digits, 0 through 9, to represent any value. The form mentioned above may be rewritten as

\begin{equation}
\sum d_{i}B^{i}
\end{equation}

where $d_{i}$ are the digits and B is the base. This is the definition of a base. For values in a base other than ten, the base will be subscribed to it. Bases will also be bolded in this text. Twelve in base \textbf{2} for example is $1100_{2}$. The first thirty-six digits that are used in this paper are in Table 1.

\section*{History}
\begin{wraptable}{l}{55mm}
  \label{digTab}
  \begin{tabular}{|c|c|@{}|c|c|}
  \hline
  Digit&Value&Digit&Value\\
  \Xhline{0.8pt}
  0&0&I&18\\1&1&J&19\\
  2&2&K&20\\3&3&L&21\\
  4&4&M&22\\5&5&N&23\\
  6&6&O&24\\7&7&P&25\\
  8&8&Q&26\\9&9&R&27\\
  A&10&S&28\\B&11&T&29\\
  C&12&U&30\\D&13&V&31\\
  E&14&W&32\\F&15&X&33\\
  G&16&Y&34\\H&17&Z&35\\
  \hline
  \end{tabular}
 \caption{First thirty-six digits and their corresponding values in base \textbf{10}.}
\end{wraptable}Numbers represented in a positional system can be traced to antiquity. Ancient Babylonians used a base sixty numeral system. The current decimal system was first developed in India around 600 A.D. At the time it was only applied to integers. Extension to the system to handle fractions wasn't until the 15th century. See Ref. \cite{knu2}, \S4.1 for more details. A. J. Kempner~\cite{kemp} gave a method to allow conversion to a base greater than 1 in 1936. A. R\'{e}nyi~\cite{reny} defined beta expansions in 1957, which is also a method to convert to positive real bases.

Bases that are less than -1 were also considered. Negative bases have the interesting property that all values can be represented without a sign~\cite{knu2}. V. Gr\"{u}nwald~\cite{grun} first discovered such bases in 1885 but was largely ignored. Kempner~\cite{kemp} suggested that such negative bases were possible. It was rediscovered again in 1957~\cite{pawl} by Pawlak and Wakulicz. These discoveries only dealt with negative bases that were integers. It wasn't until 2009 that S. Ito and T. Sadahiro came up with negative beta expansions~\cite{ito}. This method allows conversion to any real base less than -1.

Bases that are complex have also been considered, but there are few. Walter Penny discussed a number system\cite{penn} in 1965 using base \textbf{\textit{i}-1} which uses only the digits 0 and 1. Donald Knuth developed the quarter-imaginary system in 1960~\cite{knu1}. His work used \textbf{2\textit{i}} as the base and it used the digits 0, 1, 2 and 3. Such complex bases can represent every complex number without a sign~\cite{knu2}.

\section*{Conversions}
\subsection*{Positive}
The digit set for positive integer bases is defined as $\{0,\: 1, \: \ldots ,\: B-1\}$. The base itself must be two or greater. To determine the $d_{i}$ of a number N in base B where both N and B are positive integers one takes the modulo of N from B and divide N by B until N is zero. The division is rounded down to the nearest integer. N mod B for each iteration are the $d_{i}$'s. The $d_{i}$'s are sorted from last found to first. If N is negative, negate N before the algorithm and negate the result after completing the algorithm. Figure~\ref{fig:Pyth} is Python code showing this method.

For conversion to real bases greater than one a different method is needed. This method also allows conversion of real numbers to integer bases The digit set for real positive bases is defined as $\{0,\: 1,\: \ldots,\: \left \lceil B \right \rceil-1\}$. It is defined for $B>1$ and $0 \leqslant N < 1$. In the event that N is not between 0 and 1, N is divided by a suitable power of B so that it does and the algorithm is started at that power. The algorithm uses a transformation function defined as $T_B(x)=Bx-\left \lfloor Bx  \right \rfloor$. This function is called recursively such that $T_B^0 (x)=x$ , $T_B^1 (x)=T_B (x)$, $T_B^2 (x)= T_B (T_B (x))$ and so on. Therefore the digits of N in base B are found using
\begin{equation}
N = \sum_{\substack{i=-\infty\\j=0}}^{\substack{i=p-1 \\ j=\infty }} d_{i}B^{i} \textup{ with } d_{i}=\left \lfloor T_{B}^{j}(x) \right \rfloor \textup{ where } x= ^{N}/_{B^{p}} \textup{ such that } x \in \left[ 0,\: 1\right)
\end{equation}
where $i$ decreases by 1 and $j$ increases by 1 for each iteration. Should the $d_i$'s become an infinite trail of zeros, the algorithm can stop and the trailing zeros are truncated. This method is the beta expansions defined in ref.~\cite{reny}.

\subsection*{Negative}
Negative bases are also possible. The digit set for negative integer bases is $\{0,\: 1, \ldots \: ,|B|-1\}$. To convert a number N to base B where B is an integer less than 1, N is repeatedly divided by B until N becomes 0 with the digit being the modulus of B to N. The division here must be floor division.\footnote{The sign of the modulus when working with negative values is that of the divisor (the base). Not all programming languages follow this rule. See ``modulo operation'' on Wikipedia for a list of what programming languages respect this. A modulo operation that keeps the sign is equivalent to $mod(a,b)= a - (b \left \lfloor \frac{a}{b} \right \rfloor )$.} In the event that N mod B is negative, 1 is added to N and B (being negative) is subtracted from the digit. The Python code shown in Figure~\ref{fig:Pyth} also contains this method.

Conversion of real values to negative real or negative integer bases can be done using negative beta expansions~\cite{ito}. The conversion works similar to the positive real base conversion. The digit set for real negative bases is defined as $\{0,\: 1,\: \ldots,\: \lceil |B| \rceil-1\}$. It is defined for $-B < -1$ (where B is positive) and N must fall between $^{-B}/_{B+1}$ and $^{1}/_{B+1}$. These two values are defined as $l_B$ and $r_B$. In the event that N does not fall within this range, N is divided by a suitable power of B (denoted as p) so that it will fit and the algorithm is started at that power.  The transformation function is defined as $T_B (x)= -Bx-\left \lfloor -Bx-l_b \right \rfloor$. This operates similar to the transformation function for positive real bases; it is called recursively. Thus the digits are found using

\begin{equation}
N = \sum_{\substack{i=-\infty\\j=0}}^{\substack{i=p-1 \\ j=\infty }} d_{i}B^{i} \textup{ with } d_{i}=\left \lfloor -B T_{B}^{j}(x)-l_B \right \rfloor \textup{ where } x= ^{N}/_{(-B)^{p}} \textup{ such that } x \in \left[ l_B,\: r_B\right)
\end{equation}

Where $i$ decreases by 1 and $j$ increases by 1 for each iteration. In the event that one of the $d_i$’s is equal to B, 1 is subtracted from the $d_i$ and a 0 is inserted in after it. For example, if you were working with base \textbf{-8} and when you converted a number you got $0.8888888888_{-8}$; you would subtract one from each 8 and insert a zero between them leaving you with this: $0.70707070707_{-8}$.\\
\\
Conversion example: Convert 3.5 to base \textbf{-4}.
\begin{alignat*}{3}
l_b=&\, -0.8\\
r_b=&\, 0.2\\
  p=&\, 3\\
  x=& \dfrac{3.5}{(-4)^{3}}=-0.0546875\\
d_2=&\lfloor -4T_{-4}^0 (x)-l_b \rfloor=\lfloor-4x-l_b \rfloor &&=1\\
d_1=&\lfloor-4T_{-4}^1 (x)-l_b \rfloor=\lfloor-4T(x)-l_b \rfloor &&=3\\
d_0=&\lfloor-4T_{-4}^2 (x)-l_b \rfloor=\lfloor-4T(T(x))-l_b \rfloor &&=0\\
d_{-1}=&\lfloor-4T_{-4}^3 (x)-l_b \rfloor=\lfloor-4T(T(T(x)))-l_b \rfloor &&=2\\
\end{alignat*}

Because $T_{-4}^4 (x)=0$ all following digits will be zero, the process can stop and thus $3.5$ in base \textbf{-4} is $130.2_{-4}$.

\subsection*{Imaginary}
An imaginary base was first discovered by Donald Knuth in 1960\cite{knu1}. Because even powers of purely imaginary numbers can be negative numbers, an imaginary number base can depend on negative bases. This method can be extended using negative beta expansions\cite{ito} such that any purely imaginary number $C$ where $|C|>1$ can be used as a base. To do so, some modifications must be made. The digit set for an imaginary base is not that of positive or negative bases but rather $\{0,\: 1,\: \ldots,\: \lfloor |C^2| \rfloor -1\}$ with $C$ being the base. The imaginary base C has an effective negative base, $_{e}B$, that is equal to $-|C^2|$. Conversion of a complex number $A = q + r i$ where $q$ and $r$ are real numbers to an imaginary base $C$ works as follows:
\begin{enumerate}
\item Find the effective negative base, $_{e}B$, for $C$.
\item Convert $q$ to base $_{e}B$ and convert $r /\operatorname{Im}(C)$ to base $_{e}B$.
\item Adjust converted $q$ and $r$’s whole and fractional parts to be of even and odd length by leading and trailing zeros respectively.
\item Interweave $q$ with $r$ such that $q$ takes the even power positions and $r$ takes the odd power positions.
\end{enumerate}

Following these rules, base $\bm{\pi i}$ would have the effective negative base of $-\pi^2$ and use ten digits. $e^{i \frac{\pi}{6}}$ in base $\bm{\pi i}$ for instance would be $11.92771330974150459993534912112\ldots$ Additionally, an interesting relationship between positive and negative imaginary bases can be seen. For instance, let $E = 85+47i$ and $F = 85-47i$. In base \textbf{6\textit{i}} ($_{e}B = -36$) $E$ and $F$ become $\textup{10Y8D.6}_{6i}$ and $\textup{11YTD.U}_{6i}$.In base \textbf{-6\textit{i}} $E$ and $F$ become $\textup{11YTD.U}_{-6i}$ and $\textup{10Y8D.6}_{-6i}$. The sign of the base simply flips the sign of the imaginary part of the converted number.
\\
Conversion example: Convert $-5 +7i$ to base \textbf{2\textit{i}}\\
$-5$ in base \textbf{-4} is $23_{-4}$.
$7$ divided by $2$ is $3.5$, which in base $-4$ is $130.2_{-4}$.
Adding leading and trailing zeros gives ${\color{red}0023.0}_{-4}$ for the real part and ${\color{DarkGreen}0130.2}_{-4}$ for the imaginary part. Interweaving these two such that the real part takes the even power positions and the imaginary part takes the odd power positions, the result is ${\color{DarkGreen}0}{\color{red}0}{\color{DarkGreen}1}{\color{red}0}{\color{DarkGreen}3}{\color{red}2}{\color{DarkGreen}0}{\color{red}3}.{\color{DarkGreen}2}{\color{red}0}$. Thus $-5 +7i$ in base \textbf{2\textit{i}} is $103203.2_{2i}$. Converting it back to base 10 using the definition of a base to check: $1(2i)^5+ 0(2i)^4+ 3(2i)^3+ 2(2i)^2+ 0(2i)^1+ 3(2i)^0+ 2(2i)^{-1} = $
$1(32i)+ 0(16)+ 3(-8i)+ 2(-4)+ 0(2i)+ 3(1)+ 2(^{-1}/_{2i}) =$
$32i  -24i -8 +3 –i =$
$-5 +7i$.

\section*{Operations}
\subsection*{Comparisons}
While comparing complex values is not defined, one can compare the real and imaginary parts separately. In a negative base\cite{grun}\cite{eng_trans}, one compares the digits directly if they are in an even power position and oppositely in an odd power position. That is, $12_{-10} < 13_{-10}$ as $2$ and $3$ are in an even position and $2<3$ whereas $47_{-10}<27_{-10}$ as $4$ and $7$ are in an odd position and $7<4$ when compared oppositely.
Extending this to imaginary bases, real values would be compared directly in the power positions of $0,\: 4,\: 8,\:\ldots, \:4k$ (where $k$ is an integer) and compared oppositely in the positions of $2,\: 6,\: 10,\: \ldots,\: 4k+2$. For imaginary values, direct comparison would be in the $1,\: 5,\: 9,\: \ldots, \:4k+1$ positions and oppositely in the $3,\: 7,\: 11,\: \ldots, \:4k+3$ positions.

\hfill \\
\noindent Comparison example: compare the real and imaginary parts of $11873.3_{3i}$ and $10880.3_{3i}$.

The real parts are $10803_{3i}$ and $10800_{3i}$. Since the values differ in the zero position one can directly compare and thus $10803_{3i}$ is larger. The imaginary parts are $1070.3_{3i}$ and $0080.3_{3i}$. Here, the values differ in the third position and so an opposite comparison states that $80.3_{3i}$ is larger.

\subsection*{Unities}
Unless otherwise noted, the following is for imaginary bases $ni$ where $n$ is an integer and $|n|>1$. $1$ in any base is always $1$. $-1$ in any such base is $10[n^2-1]$. For example in base \textbf{3\textit{i}}: $n^2= 9$, $-1$ in this base would be $108_{3i}$. $i$ in any such bases where $n>0$ is $10.[n^2-|n|]$. $i$ in base \textbf{2\textit{i}} would be $10.2_{2i}$, in base \textbf{4\textit{i}} it would be $\textup{10.C}_{4i}$. In any such bases where $n<0$, $i$ instead becomes $0.[|n|]$. For instance, $i$ in base \textbf{9\textit{i}} would be $0.9_{9i}$. For $-i$, the rules for $+i$ above are swapped: $-i$ follows the positive integer imaginary base rule of $+i$ for negative imaginary bases and negative integer imaginary base rule of $+i$ for positive imaginary bases.

\subsection*{Arithmetic}
For addition, subtraction and multiplication they are carried out in the same manner as normal arithmetic except for a change with the carry and borrow digit rules. A carry digit removes 1 from two columns over and a borrow digit adds 1 to two columns over\cite{knu1}. Each column should not exceed the effective base. This applies to the same operations in a negative base, except they move one column over instead\cite{grun,eng_trans}. In the event that a carry digit cannot remove 1 from the column, a borrow digit is added to two columns over and the carry digit then removes 1 from the base. This can be seen in the first example below.

\hfill\break
\noindent Addition examples (all in base \textbf{3\textit{i}}, $_{e}B=-9$):\\
\hfill\break
\noindent\begin{tabu} to 0.9\textwidth{*{4}{X[c]}}
\emph{i} & \emph{ii} &\emph{iii} & \emph{iv}  \\
$\begin{array}{r}\\ 41\\ \underline{\mbox{}+  \phantom{00}61}\\ 108012\\\newline\end{array}$ &
$\begin{array}{r}\\ 132\\ \underline{\mbox{}+  11873}\\ 15\\\newline\end{array}$ &
$\begin{array}{r}\\ 0.08\\ \underline{\mbox{}+  \phantom{0}0.01}\\ 108.00\\\newline\end{array}$ & 
$\begin{array}{r}\\ 123.485\\ \underline{\mbox{} +  300.034}\\ 422.320\\\newline\end{array}$\\
\end{tabu}
%

\hfill\break
\noindent Subtraction examples (all in base \textbf{3\textit{i}}):\\
\hfill\break
\noindent\begin{tabu} to 0.9\textwidth{*{4}{X[c]}}
\emph{i} & \emph{ii} &\emph{iii} & \emph{iv}  \\
$\begin{array}{r}\\ 871\\ \underline{\mbox{}-  233}\\ 747\\\newline\end{array}$ &
$\begin{array}{r}\\ 204.000\\ \underline{\mbox{}- \phantom{00}1.104}\\ 203.005\\\newline\end{array}$ &
$\begin{array}{r}\\ 25763.0\\ \underline{\mbox{}- 126742.3}\\ 8031.6\\\newline\end{array}$ & 
$\begin{array}{r}\\ 468.782\\ \underline{\mbox{}- 551.123}\\ 10817.768\\\newline\end{array}$\\
\end{tabu}
%

\hfill\break
\noindent Multiplication examples (first three are in base \textbf{3\textit{i}}, last one is in base \textbf{4\textit{i}})\\
\hfill\break
\noindent\begin{tabu} to 0.9\textwidth{*{4}{X[c]}}
\emph{i} & \emph{ii} &\emph{iii} & \emph{iv}  \\
$\begin{array}{r}\\ 5\\ \underline{\mbox{}\times  \phantom{c0}2}\\ 10801\\\newline\end{array}$ &
$\begin{array}{r}\\ 5.0\\ \underline{\mbox{}\times \phantom{00}0.3}\\ 1080.6\\\newline\end{array}$ &
$\begin{array}{r}\\ 10432.567\\ \underline{\mbox{}\times \phantom{000}87.200}\\ 20853035\\ 147001364\phantom{0}\\ \underline{\mbox{}+ 157432732\phantom{00}}\\523204.0875\\ \newline\end{array}$ & 
$\begin{array}{r}\\ 18.68\\ \underline{\mbox{}\times \phantom{0}26.00}\\ \textup{10E4D4}\\ \underline{\mbox{}+\textup{\phantom{00}2FC\phantom{0}}}\\\textup{11FF39.4}\\ \newline\end{array}$\\
\end{tabu}
\hfill\break
%

For division, one first must multiply both the divisor and the dividend by the divisor's complex conjugate. Afterwards, the division of real and imaginary parts of the dividend can be done directly in their effective base. The real and imaginary parts of the quotient are then added back together while keeping the correct positions. This follows complex division for normal arithmetic.

If the divisor is purely real then the multiplication by the complex conjugate can be skipped. If the divisor is purely imaginary the multiplication step can be skipped but the division of the dividend's real part will result in the quotient's imaginary part and the dividend's imaginary part will result in the quotient's real part with the radix shifted one column to the left.

\hfill \break
As the effective base is a negative base, the rules for doing a division in such a base $-B$(where $B$ is positive and real) are as follows\cite{grun,eng_trans}:
\begin{enumerate}
\item The first partial dividend should be greater than the absolute value of the divisor and less than the absolute value of the divisor multiplied by $B-1$. This is the same requirement for the first partial dividend as division in a positive base.
\item If the divisor and first partial dividend are opposite in sign then the first partial quotient will be two digits long, otherwise it will be one. Once this first partial quotient is found, division may continue on as regular division. That is, each following step will have one digit.
\item Each partial remainder found must be opposite in sign of the divisor and its absolute value be smaller than the absolute value of the divisor.
\end{enumerate}

The following are two examples of of division in base \textbf{-10}:

\noindent\begin{tabu} to 0.9\textwidth{*{2}{X[c]}}
\emph{i} & \emph{ii}\\
$\begin{array}{l}\\ 
\phantom{00}\underline{\phantom{00}1512.1247}\\
28 \big) 14117\\
\phantom{00}\underline{\phantom{0}140}\\
\phantom{00000}11\\
\phantom{0000}\underline{\phantom{0}28}\\
\phantom{000000}37\\
\phantom{00000}\underline{\phantom{0}36}\\
\phantom{0000000}10\\
\phantom{000000}\underline{\phantom{0}28}\\
\phantom{00000000}20\\
\phantom{0000000}\underline{\phantom{0}36}\\
\phantom{000000000}40\\
\phantom{00000000}\underline{\phantom{0}52}\\
\phantom{0000000000}80\\
\phantom{000000000}\underline{\phantom{0}96}\\
\phantom{00000000000}4\\\newline
\end{array}$ &

$\begin{array}{l}\\
\phantom{0}\underline{\phantom{000}2261}\\
9 \big) 197349\\
\phantom{0}\underline{\phantom{0}198}\\
\phantom{000}193\\
\phantom{00}\underline{\phantom{0}198}\\
\phantom{0000}154\\
\phantom{000}\underline{\phantom{0}154}\\
\phantom{000000}09\\
\phantom{000000}\underline{\phantom{0}9}\\
\newline\end{array}$\\

\end{tabu}
\hfill\break

From this, division of $18.68$ by $\textup{10E6}$ in base \textbf{4\textit{i}} is as follows:
Complex conjugate of $\textup{10E6}$ is $\textup{10E6} \times \textup{10F} = 26$.\\
$18.68 \times 26 = \textup{11FF39.4}$\\
$\textup{10E6} \times 26 = \textup{10A04}$\\
Dividend's real part (odd positions dropped) = 1F9\\
Dividend's imaginary part (even positions dropped) = 1F3.4\\

\noindent\begin{tabu}{*{2}{X[c]}}
Real part & Imag part\\ 
$\begin{array}{l}\\ 
\phantom{000}\underline{\phantom{0000}\textup{1.C}}\\
\textup{1A4} \big) \textup{1F9}\\
\phantom{000}\underline{\phantom{0}\textup{1A4}}\\
\phantom{00000}550\\
\phantom{0000}\underline{\phantom{0}550}\\
\newline
\end{array}$ &
$\begin{array}{l}\\ 
\phantom{000}\underline{\phantom{0000}\textup{1.D}}\\
\textup{1A4} \big) \textup{1F3.4}\\
\phantom{000}\underline{\phantom{0}\textup{1A4}}\\
\phantom{00000}\textup{6F4}\\
\phantom{0000}\underline{\phantom{0}\textup{6F4}}\\
\newline
\end{array}$
\end{tabu}
\noindent$\textup{1.C}$ in the even positions $ = \textup{1.0C}$\\
$\textup{1.D}$ in the odd positions $ = \textup{10.D}$\\
$\textup{1.0C} + \textup{10.D} = \textup{11.DC}$\\

This can also be done without separating the real and imaginary parts as above and instead done entirely in line:\\
$\begin{array}{l}\\ 
\phantom{00000}\underline{\phantom{000000}\textup{11.DC}}\\
\textup{10A04} \big) \textup{11FF39.4}\\
\phantom{000000}\underline{\phantom{0}\textup{10A04}}\\
\phantom{00000000}\textup{16FF9}\\
\phantom{0000000}\underline{\phantom{0}\textup{10A04}}\\
\phantom{000000000}\textup{65F5\phantom{.}4}\\
\phantom{00000000}\underline{\phantom{0}\textup{60F0\phantom{.}4}}\\
\phantom{000000000000} 50\phantom{.}5\\
\phantom{00000000000}\underline{\phantom{0}50\phantom{.}5}\\
\newline
\end{array}$

\noindent Another example of division in base \textbf{4\textit{i}}: 11EE15FEC.168 divided by E94.\\
\noindent
$\textup{E94}^{*} = \textup{1E74}$\\
$\textup{11EE15FEC.168} \times \textup{1E74} = \textup{2143E044F7E.C}$\\
$\textup{E94} \times \textup{1E74} = \textup{104030E00}$\\
$\begin{array}{l}\\ 
\phantom{000000000}\underline{\phantom{0000000000}\textup{32A.F12}}\\
\textup{104030E00} \big) \textup{2143E044F7E.C}\\
\phantom{000000000}\underline{\phantom{0}\textup{30C070A00}}\\
\phantom{00000000000}\textup{18380A4F7}\\
\phantom{0000000000}\underline{\phantom{0}\textup{208050C00}}\\
\phantom{000000000000}\textup{8C8CA8F7E}\\
\phantom{00000000000}\underline{\phantom{0}\textup{807050C00}}\\
\phantom{00000000000000}\textup{C1C4837E\phantom{.}C}\\
\phantom{0000000000000}\underline{\phantom{0}\textup{C0A00020\phantom{.}0}}\\
\phantom{000000000000000}\textup{124835E\phantom{.}C0}\\
\phantom{00000000000000}\underline{\phantom{0}\textup{104030E\phantom{.}00}}\\
\phantom{0000000000000000}\textup{208050\phantom{.}C00}\\
\phantom{000000000000000}\underline{\phantom{0}\textup{208050\phantom{.}C00}}\\
\newline
\end{array}$
\\
An alternate method for division (which produces the inverse of the divisor) is given by Nadler~\cite{alt_div}.

\section*{Uses}
An imaginary numeral system may be useful in computing and manipulating complex numbers on computers, as previous work has stated~\cite{penn,knu1}. An alternate use for imaginary bases (and by extension all other bases) is that it can be thought of as a form of encryption. One can represent any text as a number in a base with a sufficient digit set. ``HELPIAMASTRINGOFWORDS'' can be seen as a number in base \textbf{33} or higher, base \textbf{-33} or lower, or bases $\bm{\pm5.75i}$ and larger in magnitude. This would then be converted to a different base and the converted number would be the encrypted message. In order to be able to encrypt or decrypt a message using bases, one would need to know the base the message can be read in, the base the message was converted to and the order of the digit set used.
\section*{Conclusion}
Given above are methods to convert a number into a positive, negative and imaginary base. Additionally, the four basic arithmetic operations are explained for imaginary bases. Uses for number bases are mentioned. The reader may wish to try such calculations in imaginary bases in addition to the examples given. 

\bibliographystyle{unsrt}

\begin{figure}
\begin{verbatim}
def rebase(num, base):
    "converts an integer from base ten to given integer base"

    digits = "0123456789ABCDEFGHIJKLMNOPQRSTUVWXYZ"
    values, sign = [], ""
    
    if abs(base) < 2 or abs(base) > len(digits):  # base check
        raise ValueError("invalid base")
    
    if num < 0 and base > 0: num *= -1; sign = "-"  # handle negatives
        
    while num:
        num, d = divmod(num, base)
        if d < 0: num += 1; d -= base  # handle negative base conversion
        values.append(d)

    values.reverse()
    converted = sign + "".join(digits[i] for i in values)
    return converted
\end{verbatim}
\caption{Python code to convert integers from base ten to either positive or negative integer bases}
\label{fig:Pyth}
\end{figure}

\end{document}